\documentclass[12pt]{article}
\usepackage{amsxtra,amssymb,amsthm,amsmath,latexsym}

\theoremstyle{plain}

\newtheorem{theorem}{Theorem}[section]
\numberwithin{equation}{section}

\newcommand{\refT}[1]{Theorem~\ref{T:#1}}
\newcommand{\refS}[1]{Section~\ref{S:#1}}

\def\qed{{\hfill $\Box$}}

\def\R{{\mathbb R}}

\def\doteta{\dot{\eta}}

\def\dotg{\dot{g}}
\def\doth{\dot{h}}
\def\dotu{\dot{u}}

\def\oH{{\overset{\circ}{H}}}
\def\oH1{{\overset{\circ}{H}\kern-.02in{}^1}}

\def\const{{\hbox{\,const\,}}}

\def\bee{\begin{equation*}}
\def\eee{\end{equation*}}
\def\be{\begin{equation}}
\def\ee{\end{equation}}

\begin{document}
%\begin{titlepage}
\title{   Dynamical systems method and a 
homeomorphism theorem
}

\author{
A.G. Ramm\\
 Mathematics Department, Kansas State University, \\
 Manhattan, KS 66506-2602, USA\\
ramm@math.ksu.edu\\}
%http://www.math.ksu.edu/\,$\widetilde{\ }$\,ramm}

\date{}

\maketitle\thispagestyle{empty}

\begin{abstract} \footnote{2000 Math subject classification:
6J15, 47H17, 58C15} \footnote{key words:  nonlinear
equations, homeomorphism, surjectivity, dynamical systems
method (DSM)}

Let $F$ be a nonlinear map in a real Hilbert space $H$.  
Suppose that $\sup_{u\in B(u_0,R)}$ $\|[F'(u)]^{-1}\|\leq
m(R)$, where $B(u_0,R)=\{u:\|u-u_0\|\leq R\}$, $R>0$ is
arbitrary, $u_0\in H$ is an element.  If
$\sup_{R>0}\frac{R}{m(R)}=\infty$, then $F$ is surjective.  
If $\|[F'(u)]^{-1}\|\leq a\|u\|+b$, $a\geq 0$ and $b>0$ are
constants independent of $u$, then $F$ is a homeomorphism of
$H$ onto $H$. The last result is known as an Hadamard-type
theorem, but we give a new simple proof of it based on the
DSM (dynamical systems method). \end{abstract}
%\end{titlepage}

\section{Introduction}\label{S:1}

The emphasis in this paper is on the demonstration of the
power of the DSM (dynamical system method) introduced and
applied to solving nonlinear operator equations in [3].
In this short note we give a new proof of an 
Hadamard-type
theorem on global homeomorphisms and a sufficient condition
for surjectivity of a nonlinear map in a Hilbert space. In
both cases the proof is based on the dynamical systems 
method (DSM). Although the
global homeomorphism theorem that we prove is not new, but
its proof is much shorter and simpler than the published
ones (cf [2] for example).

J.~Hadamard \cite{H} proved that a smooth map
$F:\R^n\to\R^n$ with the property $\|[F'(u)]^{-1}\|\leq b$,
$b=\const>0$, $\forall u\in\R^n$, where $F^{(j)}$ denotes
the Fr\`echet derivative, is a global homeomorphism of
$\R^n$ onto $\R^n$. This result was generalized to Hilbert
and Banach spaces under the weaker condition \be\label{e1.1}
   \|[F'(u)]^{-1}\|\leq a\|u\|+b \ee
where $a>0$ and $b>0$ are constants and $u$ is any element of the space (see \cite{OR} and references therein).  Published proofs of such a result are relatively long (cf. \cite{OR}).
In \cite{R} the DSM  (dynamical systems method) was developed as a tool for a study of nonlinear operator equations.

The aim of this paper is to apply the DSM for a proof of the following:

\begin{theorem}\label{T:1}
Assume that $F:H\to H$ is a map in a real Hilbert space and
\be\label{e1.2}
  \sup_{u\in B(u_0,R)} \|F^{(j)}(u)\|\leq M_j(R), \quad 
1\leq j\leq 2, \ee
\be\label{e1.3}
  \sup_{u\in B(u_0,R)} \|[F'(u)]^{-1}\|\leq m(R) \ee
 where $R>0$ is arbitrary and $u_0$ is an element of $H$.
 
 If 
 \be\label{e1.4}
   \sup_{R>0}\frac{R}{m(R)}=\infty \ee then $F$ is surjective. 

If \eqref{e1.1} holds, then $F$ is a global homeomorphism of $H$ onto $H$.
\end{theorem}

{\bf Remark:} Condition (1.4) is essential. For example, if 
$F(u):=e^u$,
$H=\R^1$, then equation $e^u=0$ does not have a solution, conditions (1.2) 
and (1.3) hold, but (1.4) does not hold: $m(R)=e^R$.

In \refS{2} proofs are given.

\section{Proofs}\label{S:2}
Consider the problem
\be\label{e2.1}
  \dotu=-[F'(u)]^{-1}[F(u)-f], \quad u(0)=u_0;\quad
\dotu:=\frac{du(t)}{dt}, \ee where $f\in H$ is an arbitrary given element.
From \eqref{e1.2} and \eqref{e1.3} it follows that \eqref{e2.1} has a
unique local solution. Using \eqref{e1.4} we prove that this solution is
global, i.e., exists for all $t>0$, by proving a uniform bound
$\sup_{t>0}\|u(t)\|<c$. By $c$ various positive constants are denoted.
Furthermore, we prove that $u(\infty):=\lim_{t\to\infty}u(t)$ exists, and
$F(u(\infty))=f$. Here $f\in H$ is arbitrary, so $F$ is surjective. The
above scheme is the dynamical systems method (DSM). Let us give the
details.

Denote $\|F(u(t))-f\|:=g(t)$. Then, by \eqref{e2.1}
$g \dotg=-g^2$, so $g(t)\leq g(0)e^{-t}$, and
\be\label{e2.2}
  \|\dotu\|\leq m(R)g(0)e^{-t}.\ee If the solution $u(t)$
does not leave the ball $B(u_0,R)$ for all times, then
$u(t)$ exists for all $t>0$. Integrating \eqref{e2.2} yields
$\|u(t)-u(0)\|\leq m(R)g(0)$. If there is an $R>0$ such that
\be\label{e2.3}
  m(R)g(0)\leq R,\ee then $u(t)\in B(u_0,R)$ $\forall t>0$,
so $u(t)$ is the global solution to \eqref {e2.1}. Condition
\eqref{e1.4} guarantees that for any fixed $u_0$ there is an
$R>0$ such that \eqref{e2.3} holds. For this $R$ one has
$u(t)\in B(u_0,R)$ $\forall t>0$, there exists $u(\infty)$,
the following estimate holds:
\be\label{e2.4}
  \|u(t)-u(\infty)\|\leq m(R)g(0)e^{-t},\ee
and, passing to the limit $t\to\infty$ in 
\eqref{e2.1}  yields $F(u(\infty))=f$.

{\it This proves the surjectivity of $F$.}

If \eqref{e1.1} holds, then \eqref{e2.2} is replaced by
\be\label{e2.5}
  \|\dotu\|\leq(a\|u(t)\|+b)g(0)e^{-t}. \ee Let $h(t):=\|u(t)\|$. Then
\eqref{e2.5} yields 
$$\doth\leq(h+p)ag_0e^{-t},\quad p:=\frac{b}{a}>0,
$$ 
so
$$\sup_{t\geq 0}h(t)\leq c_1:=(\|u(0)\|+p)e^{ag_0}-p,
$$ 
and 
\be\label{e2.6}
  \|\dotu\|\leq c_2 e^{-t}, \quad c_2:=(ac_1+b)g(0).\ee Thus
$u(t)\in B(u_0,c_2)$. It is well known and easy to prove
that condition \eqref{e1.3} implies that $F$ is a local
homeomorphism, i.e., $F$ maps a neighborhood of any point
$u\in H$ homeomorphically onto a neighborhood of a point
$F(u)$. From \eqref{e2.1} and the estimate $\|\dotu\|\leq
c_2e^{-t}$ we conclude as above that $F$ is surjective.
Therefore, in order to prove that $F$ is a global
homeomorphism of $H$ onto $H$ it is sufficient to prove that
$F(u)=F(v)$ implies $u=v$.

{\it The idea of our proof} is to consider the path
$(1-s)u_0+sv:=w(s):=u(0,w(s)):=u(0,s)$ from $u_0$ to $v$ and 
to construct the
solution $u(t,w(s)):=u(t,s)$ to \eqref{e2.1} with the
initial data $w(s)$ replacing $u_0$, and then to show that
$u(\infty,s)=u$ for each $s$ and to conclude that
$v=u(\infty,1)=u$. In this step we use the assumption(1.3) 
which implies that $F$ is a local
homeomorphism: if $F(u(\infty,s))=F(u(\infty,s+\sigma))=f$
and $||u(\infty,s)-u(\infty,s+\sigma)||$ is sufficiently
small, then $u(\infty,s)=u(\infty,s+\sigma)$ provided that $F$ is a
local homeomorphism.

Let us give the details.

If $s=0$ then we have $u(\infty,0)=u$. If $\sigma$ is small, then 
\be\label{e2.7}
  \sup_{t\geq 0} \|u(t,s+\sigma)-u(t,s)\|\leq c
\|u(0,s+\sigma)-u(0,s)\|, \ee where $c>0$ does not depend on
$s$, $r$, and $t$, and $\sigma>0$ does not depend on $s$.  
We prove \eqref{e2.7} below. If \eqref{e2.7} holds, then
$\|u(\infty,s+\sigma)-u(\infty,s)\|$ is arbitrarily small if
$\|u(0,s+\sigma)-u(0,s)\|:=\delta$ is sufficiently small.
Since $F(u(\infty,s+\sigma))=F(u(\infty,s))  =f$, and since
$F$ is a local homeomorphism, it follows that
$u(\infty,s+\sigma)=u(\infty,s)$. Since $u(\infty,0)=u$, and
since $\sigma>0$ does not depend on $s$, then in finite number
of steps one gets to the point $s+\sigma=1$ and concludes
that $u=u(\infty,s)=u(\infty,1)=v$, $0\leq s\leq 1$. Thus,
to complete the proof we have to check \eqref{e2.7}. Denote
$x(t):=u(t,s+\sigma)-u(t,s):=z-y$, and $\|x(t)\|:=\eta(t)$.
Then, using \eqref{e2.6} and \eqref{e1.2}, one gets:
\be\label{e2.8}
  \begin{aligned}
  \eta\doteta&=-([F'(z)]^{-1}(F(z)-f)
-[F'(y)]^{-1}(F(y)-f),x(t))\\
                & =-(([F'(z)]^{-1}-[F'(y)]^{-1}) (F(z)-f),x)
-([F'(y)]^{-1}(F(z)-F(y)),x)\\
                & \quad \leq c e^{-t}\eta^2-\eta^2+c\eta^3,
\end{aligned}\ee where $c>0$ is a constant and we have used
the formulas: $$ F(z)-F(y)=F'(y)(z-y)+K,\quad
\|K\|\leq\frac{M_2}{2}\|z-y\|^2. $$ Since $\eta\geq 0$, one
gets from \eqref{e2.8} the inequality \be\label{e2.9}
  \doteta\leq -\eta+c\eta^2+ce^{-t}\eta,\quad \eta(0)=\delta. \ee
Let    $\eta=he^{-t}$. Then
\be\label{e2.10}
 \doth\leq ce^{-t}(h^2+h),\quad h(0)=\delta.\ee
From \eqref{e2.10} one gets $h(t)\leq c_3\delta$, 
so 
$$\eta(t)\leq c_3 e^{-t}\delta.$$
 This implies \eqref{e2.7}.
 
 \refT{1} is proved. \qed

\end{document}